\documentclass{amsart}

\usepackage{cite}

\numberwithin{equation}{section}

\newcommand{\abs}[1]{\lvert#1\rvert}

\begin{document}

\title[Elliptic hypergeometric series]
{Summation formulae for elliptic hypergeometric series}

\author{S.~Ole Warnaar}\thanks{Work supported by the Australian
Research Council}
\address{Department of Mathematics and Statistics,
The University of Melbourne, Vic 3010, Australia}
\email{warnaar@ms.unimelb.edu.au}

\keywords{basic hypergeometric series, elliptic hypergeometric series,
Bailey lemma}

\subjclass[2000]{33D15, 33E05}

\begin{abstract}
Several new identities for elliptic hypergeometric series are
proved. Remarkably, some of these are elliptic analogues of
identities for basic hypergeometric series that are balanced but
not very-well-poised.
\end{abstract}

\maketitle

\section{Introduction}
Recently there has been much interest in elliptic hypergeometric
series
\cite{vDS00,vDS01a,vDS01b,vDS02,FT97,KN03,KNR03,Rosengren01,
Rosengren02,RS03,Spiridonov01,
Spiridonov02a,Spiridonov02,Spiridonov03,SZ99,W02,W03}. The
simplest examples of such series are of the type
\begin{multline}\label{Vr}
{_{r+1}V_r}(a_1;a_6,\dots,a_{r+1};q,p) \\=
\sum_{k=0}^{\infty}\frac{\theta(a_1 q^{2k};p)}{\theta(a_1;p)}
\frac{(a_1,a_6,\dots,a_{r+1};q,p)_k}
{(q,a_1q/a_6,\dots,a_1q/a_{r+1};q,p)_k}\, q^k,
\end{multline}
where $\theta(a;p)$ is a theta function
\begin{equation*}
\theta(a;p)=\prod_{i=0}^{\infty}(1-ap^i)(1-p^{i+1}/a),
\qquad 0<\abs{p}<1,
\end{equation*}
and $(a;q,p)_n$ is the elliptic analogue of the $q$-shifted factorial
\begin{equation*}
(a;q,p)_n=\prod_{j=0}^{n-1}\theta(aq^j;p).
\end{equation*}
As usual,
\begin{equation*}
(a_1,\dots,a_k;q,p)_n=(a_1;q,p)_n\dots(a_k;q,p)_n.
\end{equation*}

For reasons of convergence one must impose that one of the parameters
$a_i$ is of the form $q^{-n}$ so that the above series terminates.
Furthermore, to obtain non-trivial results, $r$ must be odd and
\begin{equation*}
a_6\cdots a_{r+1}q=(a_1 q)^{(r-5)/2}.
\end{equation*}

For ordinary as well as basic hypergeometric series a vast number
of summation identities are known, see e.g., \cite{GR90,Slater66}.
Unfortunately, most of these do not appear to have an elliptic
analogue and to the best of my knowledge the only two summation
identities for series of the type \eqref{Vr} known to date are the
elliptic Jackson sum of Frenkel and Turaev 
\cite[Theorem 5.5.2]{FT97}
\begin{equation}\label{V109}
{_{10}V_9}(a;b,c,d,e,q^{-n};q,p)=
\frac{(aq,aq/bc,aq/bd,aq/cd;q,p)_n}{(aq/b,aq/c,aq/d,aq/bcd;q,p)_n},
\end{equation}
for $bcde=a^2 q^{n+1}$, and the identity \cite[Theorem 4.1]{W02}
\begin{multline*}
{_{2r+8}V_{2r+7}}(ab;c,ab/c,bq,bq^2,\dots,bq^r,aq^n,
aq^{n+1},\dots,aq^{n+r-1},q^{-rn};q^r,p) \\
=\frac{(a/c,c/b;q,p)_n} {(cq^r,abq^r/c;q^r,p)_n}
\frac{(q^r,abq^r;q^r,p)_n} {(a,1/b;q,p)_n}.
\end{multline*}

In a recent paper \cite{W03} I stated without proof
that
\begin{multline}\label{biba}
\sum_{k=0}^n \frac{\theta(a^2q^{4k};p^2)}{\theta(a^2;p^2)}
\frac{(a^2,b/q;q^2,p^2)_k}{(q^2,a^2 q^3/b;q^2,p^2)_k}
\frac{(a q^n/b,q^{-n};q,p)_k}{(bq^{1-n},aq^{n+1};q,p)_k}\,q^{2k} \\
=\frac{\theta(-aq^{2n}/b;p)}{\theta(-a/b;p)}
\frac{(-a/b,aq;q,p)_n}{(-q,1/b;q,p)_n}
\frac{(1/bq;q^2,p^2)_n}{(a^2q^3/b;q^2,p^2)_n}\,q^n.
\end{multline}
When $p$ tends to zero this simplifies to a bibasic summation of
Nassrallah and Rahman \cite[Corollary 4]{NR81} (see also
\cite[Equation (3.10.8)]{GR90}). Initially I was only able to find
a rather unpleasant inductive proof, but an e-mail exchange with
Vyacheslav Spiridonov prompted me to try again to find a more
constructive derivation of \eqref{biba}. In this paper I will give
such a proof. Interestingly, it depends crucially on the new
elliptic identity
\begin{multline}\label{new1}
{_{12}V_{11}}(ab;b,bq,b/p,bqp,aq^2/b,a^2q^{2n},q^{-2n};q^2,p^2) \\
=\frac{\theta(a;p)}{\theta(aq^{2n};p)}
\frac{(-q,aq/b;q,p)_n}{(a,-b;q,p)_n}
\frac{(abq^2;q^2,p^2)_n}{(a/b;q^2,p^2)_n}\,q^{-n},
\end{multline}
which provides a third example of a summable ${_{r+1}V_r}$ series.

The quasi-periodicity of the theta functions
\begin{equation}\label{quasi1}
\theta(a;p)=-a\, \theta(ap;p)
\end{equation}
yields
\begin{equation}\label{quasi2}
(a;q,p)_n=(-a)^n q^{\binom{n}{2}}(ap;q,p)_n.
\end{equation}
Morover, from
\begin{equation}\label{lim1}
\lim_{p\to 0} \theta(ap;p^2)=1
\end{equation}
it follows that
\begin{equation}\label{lim2}
\lim_{p\to 0} (ap;q,p^2)_n=1.
\end{equation}
Hence
\begin{equation*}
\lim_{p\to 0}\frac{(b/p;q^2,p^2)_k}{(aq/p;q^2,p^2)_k}
=\Bigl(\frac{b}{aq}\Bigr)^k
\lim_{p\to 0}\frac{(bp;q^2,p^2)_k}{(aqp;q^2,p^2)_k}
=\Bigl(\frac{b}{aq}\Bigr)^k.
\end{equation*}
Using standard notation for basic hypergeometric series
\cite{GR90} it thus follows that in the $p\to 0$ limit \eqref{new1} 
becomes
\begin{multline*}
{_8W_7}(ab;b,bq,aq^2/b,a^2q^{2n},q^{-2n};q^2,bq/a) \\
=\frac{1-a}{1-aq^{2n}} \frac{(-q,aq/b;q)_n}{(a,-b;q)_n}
\frac{(abq^2;q^2)_n}{(a/b;q^2)_n}\,q^{-n}.
\end{multline*}
Using Watson's $_8\phi_7$ transformation \cite[Equation (III.18)]{GR90}
this may be also put as
\begin{equation}\label{BW}
{_4\phi_3}\biggl[\genfrac{}{}{0pt}{}{b,bq,a^2q^{2n},q^{-2n}}
{b^2,aq,aq^2};q^2,q^2\biggr]= \frac{1-a}{1-aq^{2n}}
\frac{(-q,aq/b;q)_n}{(a,-b;q)_n}\,b^n,
\end{equation}
an identity discovered recently in \cite{BW03}.

Given \eqref{new1} the proof of \eqref{biba} is routine, but
proving \eqref{new1} is unexpectetly difficult since its 
constructive proof requires \eqref{biba}!
In the next section I will therefore give a rather non-standard proof of
\eqref{new1} by specializing a recent elliptic transformation formula 
of Spiridonov in a singular point. 
The bonus of this proof is that it immediately suggests the 
following companion to \eqref{new1}
\begin{multline}\label{new2}
{_{12}V_{11}}(ab;b,-b,bp,-b/p,aq/b,a^2q^{n+1},q^{-n};q,p^2) \\
=\chi(n\textup{ is even})
\frac{(q,a^2q^2/b^2;q^2,p^2)_{n/2}}{(a^2q^2,b^2q;q^2,p^2)_{n/2}}
\frac{(abq;q,p^2)_n}{(aq/b;q,p^2)_n},
\end{multline}
with $\chi(\text{true})=1$ and $\chi(\text{false})=0$.
This is the fourth example of a ${_{r+1}V_r}$ that can be summed.
In the limit when $p$ tends to zero \eqref{new2} simplifies to
\begin{multline*}
{_8W_7}(ab;b,-b,aq/b,a^2q^{n+1},q^{-n};q,-b/a) \\
=\chi(n\textup{ is even})
\frac{(q,a^2q^2/b^2;q^2)_{n/2}}{(a^2q^2,b^2q;q^2)_{n/2}}
\frac{(abq;q)_n}{(aq/b;q)_n}.
\end{multline*}
By Watson's $_8\phi_7$ transformation this can be further reduced
to Andrews' terminating $q$-analogue of
Watson's $_3F_2$ sum \cite[Theorem 1]{Andrews76}
(see also \cite[Equation (II.17)]{GR90})
\begin{equation}\label{A}
{_4\phi_3}\biggl[\genfrac{}{}{0pt}{}{b,-b,a^2q^{n+1},q^{-n}}
{b^2,aq,-aq};q,q\biggr]= \chi(n\textup{ is even})
\frac{(q,a^2q^2/b^2;q^2,p)_{n/2}}{(a^2q^2,b^2q;q^2,p)_{n/2}}\,b^n.
\end{equation}

The identities \eqref{new1} and \eqref{new2} together with
Watson's transformation imply the $_4\phi_3$ sums \eqref{BW}
and \eqref{A}. 
It is however also possible to rewrite \eqref{new1} and \eqref{new2}
as two elliptic summations that yield \eqref{BW}
and \eqref{A} when $p$ tends to zero without
an appeal to Watson's transformation. 
Making the substitution $a\to ap$ in \eqref{new1} and using the
quasi-periodicities \eqref{quasi1} and \eqref{quasi2} yields
\begin{multline}\label{linp}
{_{12}V_{11}}(abp;b,bq,bp,bqp,aq^2p/b,a^2q^{2n},q^{-2n};q^2,p^2) \\
=\frac{\theta(a;p)}{\theta(aq^{2n};p)}
\frac{(-q,aq/b;q,p)_n}{(a,-b;q,p)_n}
\frac{(abq^2p;q^2,p^2)_n}{(ap/b;q^2,p^2)_n}\, b^n.
\end{multline}
By \eqref{lim1} and \eqref{lim2} the $p\to 0$ limit
breaks the very-well-poisedness, resulting in \eqref{BW}. 
In much the same way, replacing $a\to ap$ in \eqref{new2} and
using \eqref{quasi1} and \eqref{quasi2} yields
\begin{multline}\label{linp2}
{_{12}V_{11}}(abp;b,-b,bp,-bp,aqp/b,a^2q^{n+1},q^{-n};q,p^2) \\
=\chi(n\textup{ is even})
\frac{(q,a^2q^2/b^2;q^2,p^2)_{n/2}}{(a^2q^2,b^2q;q^2,p^2)_{n/2}}
\frac{(abqp;q,p^2)_n}{(aqp/b;q,p^2)_n}\, b^n.
\end{multline}
When $p$ tend to $0$ this reduces to \eqref{A}.

The results \eqref{linp} and \eqref{linp2} show that, potentially, 
many more identities for series that are balanced but not very-well poised 
may have an elliptic analogue. Indeed, after showing him \eqref{linp} and 
\eqref{linp2}, Michael Schlosser observed that making
the simultaneous variable changes $\{a,d,e,p\}\to\{ap,aqp/d,ep,p^2\}$ in
\eqref{V109} gives
\begin{equation*}
{_{10}V_9}(ap;b,c,aqp/d,ep,q^{-n};q,p^2)=
\frac{(aqp,aqp/bc,d/b,d/c;q,p^2)_n}
{(aqp/b,aqp/c,d,d/bc;q,p^2)_n},
\end{equation*}
for $bce=adq^n$.
In the $p$ to $0$ limit this results in the $q$-Pfaff--Saalsch\"utz sum
\cite[Equation (II.12)]{GR90}
\begin{equation*}
{_3\phi_2}\biggl[\genfrac{}{}{0pt}{}{b,c,q^{-n}}{d,bcq^{1-n}/d};q,q\biggl]
=\frac{(d/b,d/c;q)_n}{(d,d/bc;q)_n}.
\end{equation*}

Probably the most important balanced summation not yet treated
is Andrews' terminating $q$-analogue of Wipple's $_3F_2$ sum
\cite[Theorem 2]{Andrews76} (see also \cite[Equation (II.19)]{GR90})
\begin{equation}\label{qWhipple}
{_4\phi_3}\biggl[\genfrac{}{}{0pt}{}{b,-b,q^{n+1},q^{-n}}
{-q,c,b^2q/c};q,q\biggr]
=\frac{(c/b^2;q)_n}{(c;q)_n}\frac{(cq^{-n};q^2)_n}{(cq^{-n}/b^2;q^2)_n}.
\end{equation}
To obtain its elliptic analogue I will first prove the new identity
\begin{multline}\label{new3}
{_{12}V_{11}}(b;-b,bp,-b/p,c/b,bq/c,q^{n+1},q^{-n};q,p^2) \\
=\frac{(bq,c/b^2;q,p^2)_n}{(q/b,c;q,p^2)_n}
\frac{(cq^{-n};q^2,p^2)_n}{(cq^{-n}/b^2;q^2,p^2)_n}\,(-1/b)^n.
\end{multline}
Replacing $b\to bp$ and using \eqref{quasi1} and \eqref{quasi2}
this implies the identity
\begin{multline*}
{_{12}V_{11}}(bp;b,-b,-bp,cp/b,bpq/c,q^{n+1},q^{-n};q,p^2) \\
=\frac{(bqp,c/b^2;q,p^2)_n}{(qp/b,c;q,p^2)_n}
\frac{(cq^{-n};q^2,p^2)_n}{(cq^{-n}/b^2;q^2,p^2)_n},
\end{multline*}
which simplifies to \eqref{qWhipple} when $p$ tends to $0$
thanks to \eqref{lim1} and \eqref{lim2}.

\section{Proofs of \eqref{biba}, \eqref{new1}, \eqref{new2}
and \eqref{new3}}
First I will give a proof of \eqref{biba} assuming \eqref{new1},
and a proof of \eqref{new1} assuming \eqref{biba}. Then I will
give a different proof of \eqref{new1} based on the transformation
\eqref{Sp} below.

\begin{proof}[Proof of \eqref{biba} based on \eqref{new1}]
When $cd=aq$ equation \eqref{V109} simplifies to
\begin{equation}\label{V87}
{_8V_7}(a;b,aq^n/b,q^{-n};q,p)=\delta_{n,0},
\end{equation}
with $\delta_{n,m}=\chi(n=m)$. Making the simultaneous replacements
\begin{equation*}
\{a,b,n,q,p\}\to \{a^2,b/q,r,q^2,p^2\},
\end{equation*}
then multiplying both sides by
\begin{equation*}
\frac{\theta(a^2q^{4r+1}/b;p^2)}{\theta(a^2q/b;p^2)}
\frac{(-aq;q,p)_{2r}}{(-aq/b;q,p)_{2r}}
\frac{(a^2q/b,q/b,a^2q^{2n}/b^2,q^{-2n};q^2,p^2)_r}
{(q^2,a^2q^2,bq^{3-2n},a^2q^{2n+3}/b;q^2,p^2)_r}\,(bq^2)^r
\end{equation*}
and finally summing $r$ from $0$ to $n$ yields
\begin{multline*}
\sum_{r=0}^n \frac{\theta(a^2q^{4r+1}/b;p^2)}{\theta(a^2q/b;p^2)}
\frac{(-aq;q,p)_{2r}}{(-aq/b;q,p)_{2r}}
\frac{(a^2q/b,q/b,a^2q^{2n}/b^2,q^{-2n};q^2,p^2)_r}
{(q^2,a^2q^2,bq^{3-2n},a^2q^{2n+3}/b;q^2,p^2)_r}\,
(bq^2)^r \\
\times {_8V_7}(a^2;b/q,a^2q^{2r+1}/b,q^{-2r};q^2,p^2)=1.
\end{multline*}
Interchanging the order of summation and using the
identity
\begin{equation}\label{id}
\frac{(a;q,p)_{2n}}{(b;q,p)_{2n}}=
\frac{(a,aq,a/p,aqp;q^2,p^2)_n}
{(b,bq,b/p,bqp;q^2,p^2)_n}
\Bigl(\frac{b}{a}\Bigr)^n
\end{equation}
this becomes
\begin{multline*}
\sum_{s=0}^n \frac{(-aq,q,p)_{2s}}{(-aq/b;q,p)_{2s}}
\frac{(a^2q^3/b;q^2,p^2)_{2s}}{(a^2;q^2,p^2)_{2s}}
\frac{(a^2,b/q,a^2q^{2n}/b^2,q^{-2n};q^2,p^2)_s}
{(q^2,a^2q^3/b,bq^{3-2n},a^2q^{2n+3}/b;q^2,p^2)_s}\, q^{3s} \\
\times {_{12}V_{11}}(a^2q^{4s+1}/b;-aq^{2s+1},-aq^{2s+2},
-aq^{2s+1}/p,-aq^{2s+2}p, \\
q/b,a^2q^{2n+2s}/b^2,q^{2s-2n};q^2,p^2)=1.
\end{multline*}
Summing the ${_{12}V_{11}}$ series by \eqref{new1} and making
some simplifications completes the proof.
\end{proof}

\begin{proof}[Proof of \eqref{new1} based on \eqref{biba}]
Replacing 
\begin{equation*}
\{a,b,n,q,p\}\to \{a,aq^2/b^2,r,q^2,p^2\}
\end{equation*} in \eqref{V87}, multiplying both sides by
\begin{equation*}
\frac{\theta(b^2q^{4r-2};p^2)}{\theta(b^2/q^2;p^2)}
\frac{(b^2/q^2,b^2/aq^2;q^2,p^2)_r}{(q^2,aq^2;q^2,p^2)_r}
\frac{(-aq^n/b,q^{-n};q,p)_r}
{(b^2q^{-n}/a,-bq^n;q,p)_r}\, q^{2r}
\end{equation*}
and summing $r$ from $0$ to $n$ yields
\begin{multline*}
\sum_{r=0}^n
\frac{\theta(b^2q^{4r-2};p^2)}{\theta(b^2/q^2;p^2)}
\frac{(b^2/q^2,b^2/aq^2;q^2,p^2)_r}{(q^2,aq^2;q^2,p^2)_r}
\frac{(-aq^n/b,q^{-n};q,p)_r}
{(b^2q^{-n}/a,-bq^n;q,p)_r}\, q^{2r} \\
\times {_8V_7}(a;aq^2/b^2,b^2q^{2r-2},q^{-2r};q^2,p^2)=1.
\end{multline*}
A change in the order of summation leads to
\begin{multline*}
\sum_{s=0}^n
\frac{\theta(b^2q^{4s-2};p^2)}{\theta(b^2/q^2;p^2)}
\frac{(b^2/q^2;q^2,p^2)_{2s}}{(a;q^2,p^2)_{2s}}
\frac{(a,aq^2/b^2;q^2,p^2)_s}{(q^2,b^2;q^2,p^2)_s}
\frac{(-aq^n/b,q^{-n};q,p)_s}{(b^2q^{-n}/a,-bq^n;q,p)_s} \\
\times
\Bigl(\frac{b^2}{a}\Bigr)^s  \:
\sum_{r=0}^{n-s}
\frac{\theta(b^2q^{4r+4s-2};p^2)}{\theta(b^2q^{4s-2};p^2)}
\frac{(b^2q^{4s-2},b^2/aq^2;q^2,p^2)_r}{(q^2,aq^{4s+2};q^2,p^2)_r} \\
\times
\frac{(-aq^{n+s}/b,q^{s-n};q,p)_r}
{(b^2q^{s-n}/a,-bq^{n+s};q,p)_r}\,q^{2r}=1.
\end{multline*}
The sum over $r$ can be performed by \eqref{biba} giving
\begin{multline*}
\sum_{s=0}^n
\frac{\theta(aq^{4s};p^2)}{\theta(a;p^2)}
\frac{(b;q,p)_{2s}}{(aq/b;q,p)_{2s}}
\frac{(a,aq^2/b^2,a^2q^{2n}/b^2,q^{-2n};q^2,p^2)_s}
{(q^2,b^2,b^2q^{2-2n}/a,aq^{2n+2};q^2,p^2)_s}
\Bigl(\frac{b^2q}{a}\Bigr)^s \\
=q^{-n} \frac{\theta(a/b;p)}{\theta(aq^{2n}/b;p)}
\frac{(-q,aq/b^2;q,p)_n}{(a/b,-b;q,p)_n}
\frac{(aq^2;q^2,p^2)_n}{(a/b^2;q^2,p^2)_n}.
\end{multline*}
Once more using \eqref{id} and replacing $a$ by $ab$ completes the proof.
\end{proof}

\begin{proof}[Proof of \eqref{new1}]
To give a proof of \eqref{new1} that does not rely on \eqref{biba}
I need the following transformation formula of Spiridonov
\cite[Theorem 5.1]{Spiridonov02} (see also \cite[Theorem 4.1]{W02}):
\begin{align}\label{Sp}
{_{14}V_{13}}&(a;a^2q/m,b^{1/2},-b^{1/2},c^{1/2},-c^{1/2},
k^{1/2}q^n,-k^{1/2}q^n,q^{-n},-q^{-n};q,p) \\
&=\frac{(a^2q^2,k/m,mq^2/b,mq^2/c;q^2,p^2)_n}
{(mq^2,k/a^2,a^2q^2/b,a^2q^2/c;q^2,p^2)_n} \notag \\
& \quad \qquad \times {_{14}V_{13}}(m;a^2q^2/m,d,dq,d/p,dqp,b,c,
kq^{2n},q^{-2n};q^2,p^2), \notag
\end{align}
for $m=bck/a^2q^2$ and $d=-m/a$. When $p$ tends to $0$ this becomes
\begin{align}\label{NR414}
{_{12}W_{11}}&(a;a^2q/m,b^{1/2},-b^{1/2},c^{1/2},-c^{1/2},
k^{1/2}q^n,-k^{1/2}q^n,q^{-n},-q^{-n};q,q) \\
&\quad =\frac{(a^2q^2,k/m,mq^2/b,mq^2/c;q^2)_n}
{(mq^2,k/a^2,a^2q^2/b,a^2q^2/c;q^2)_n} \notag \\
& \qquad \qquad \qquad \times
{_{10}W_9}(m;a^2q^2/m,d,dq,b,c,kq^{2n},q^{-2n};q^2,mq/a^2)
\notag
\end{align}
which is equivalent to a bibasic transformation of
Nassrallah and Rahman \cite[Equation (4.14)]{NR81}
(see also \cite[Equation (3.10.15)]{GR90}).
In the above representation \eqref{NR414} has been rediscovered
very recently in \cite[Equation (4.9)]{AB02}.

To now prove \eqref{new1} I observe that the
$_{14}V_{13}$ series on the left side of \eqref{Sp} as well as the
prefactor on the right side of \eqref{Sp} are singular for $k=a^2$. 
Multiplying both sides by
$(k/a^2;q^2,p^2)_n$ and observing that for $0\leq r\leq n$
\begin{equation*}
\lim_{k\to a^2} \frac{(k/a^2;q^2,p^2)_n}{(a^2q^{2-2n}/k;q^2,p^2)_r}
=(-1)^n q^{n^2-n}\delta_{n,r},
\end{equation*}
it follows that in the limit when $k$ tends to $a^2$ only 
the term with $r=n$ survives in the sum on the left 
(with $r$ being the summation index of the $_{14}V_{13}$ series).
As a result
\begin{multline*}
{_{12}V_{11}}(m;a^2q^2/m,d,dq,d/p,dqp,
a^2q^{2n},q^{-2n};q^2,p^2) \\
=q^{-n}\frac{\theta(-a;p)}{\theta(-aq^{2n};p)}
\frac{(-q,a^2q/m;q,p)_n}{(-a,m/a;q,p)_n}
\frac{(mq^2;q^2,p^2)_n}{(a^2/m;q^2,p^2)_n},
\end{multline*}
with $m=bc/q^2$ and $d=-m/a$.
Since the only dependence on $b$ and $c$ is through the definition
of $m$, the equation $m=bc/q^2$ is superfluous, and
the above is true with $a$ and $m$ arbitrary indeterminates.
Making the simultaneous changes $m\to ab$ and $a\to -a$ yields
\eqref{new1}.
\end{proof}

\begin{proof}[Proof of \eqref{new2}]
As mentioned in the introduction, the above proof of \eqref{new1}
immediately suggests \eqref{new2} by virtue of the fact that
\eqref{Sp} has the companion \cite[Theorem 4.2]{W02}
\begin{align}\label{Wa}
{_{14}V_{13}}&(a;a^2/m^2,b,bq,c,cq,kq^n,kq^{n+1},q^{-n},q^{1-n};q^2,p) \\
&=\frac{(aq,k/m,mq/b,mq/c;q,p)_n}{(mq,k/a,aq/b,aq/c;q,p)_n}  \notag \\
&\qquad\quad \times
{_{14}V_{13}}(m;a/m,d,-d,
dp^{1/2},-d/p^{1/2},b,c,kq^n,q^{-n};q,p), \notag
\end{align}
for $m=bck/aq$ and $d=m(q/a)^{1/2}$.
In the $p\to 0$ limit this gives 
\begin{multline}\label{RV78}
{_{12}W_{11}}(a;a^2/m^2,b,bq,c,cq,kq^n,kq^{n+1},q^{-n},q^{1-n};q^2,q^2) \\
=\frac{(aq,k/m,mq/b,mq/c;q)_n}{(mq,k/a,aq/b,aq/c;q)_n} \,
{_{10}W_9}(m;a/m,d,-d,b,c,kq^n,q^{-n};q,-mq/a)
\end{multline}
due to Rahman and Verma \cite[Equation (7.8)]{RV93}
(see also \cite[Equation (3.13)]{AB02}).

This time the singularity to be exploited occurs for $k=a$.
Multiplying both sides of \eqref{Wa} by
$(k/a;q,p)_n$ and observing that for $0\leq 2r\leq n$
\begin{equation*}
\lim_{k\to a} \frac{(k/a;q,p)_n}{(aq^{1-n}/k;q^2,p^2)_{2r}}
=q^{\binom{n}{2}}\delta_{n,2r},
\end{equation*}
it follows that in the $k\to a$ limit
only the term with $2r=n$ survives
in the sum on the left (with $r$ being the summation index of the
$_{14}V_{13}$ series). Hence
\begin{multline*}
{_{12}V_{11}}(m;a/m,d,-d,dp^{1/2},-d/p^{1/2},aq^n,q^{-n};q,p) \\
=\chi(n\text{ even})
\frac{(a,a^2/m^2;q^2,p)_{n/2}}{(q^2,m^2q^2/a;q^2,p)_{n/2}}
\frac{(q,mq;q,p)_n}{(a,a/m;q,p)_n}
\end{multline*}
with $m=bc/q$ and $d=m(q/a)^{1/2}$.
Again the dependence on $b$ and $c$ is only through the definition
of $m$, so that the above is true for arbitrary $a$ and $m$.
Making the simultaneous changes $m\to ab$, $a\to a^2q$ and $p\to p^2$
yields \eqref{new2}.
\end{proof}

\begin{proof}[Proof of \eqref{new3}]
Making the simultaneous substitutions
\begin{equation*}
\{a,b,c,d,e,f,g,p\}\to
\{b,c/b,bq/c,q^{n+1},-b,bp,-b/p,p^2\}
\end{equation*}
in the elliptic analogue of Bailey's $_{10}\phi_9$
transformation \cite[Theorem 5.5.1]{FT97}
\begin{multline*}
{_{12}V_{11}}(a;b,c,d,e,f,g,q^{-n};q,p) \\
=\frac{(aq,aq/ef,aq/fg,aq/eg;q,p)_n}
{(aq/e,aq/f,aq/g,aq/efg;q,p)_n}\:
{_{12}V_{11}}(\lambda;\lambda b/a,\lambda c/a,\lambda d/a,e,f,g,q^{-n};q,p)
\end{multline*}
for $bcdefg=a^3 q^{n+2}$ and $\lambda=a^2q/bcd$, \eqref{new3}
can be transformed into
\begin{multline}\label{V1211}
{_{12}V_{11}}(b^2q^{-n-1};b,-b,bp,-b/p,cq^{-n-1},b^2q^{-n}/c,q^{-n};q,p^2) \\
=\frac{(q/b^2,c/b^2;q,p^2)_n}{(q,c;q,p^2)_n}
\frac{(q^2,cq^{-n};q^2,p^2)_n}{(q^2/b^2,cq^{-n}/b^2;q^2,p^2)_n}.
\end{multline}
Here the right-hand side has been simplified using
\begin{equation*}
\frac{(a,-a,a/p,-ap;q,p^2)_n}{(b,-b,bp,-b/p;q,p^2)_n}
=\frac{(a^2;q^2,p^2)_n}{(b^2;q^2,p^2)_n}
\Bigl(-\frac{a}{b}\Bigr)^n
\end{equation*}
with $a\to q$ and $b\to q/b$.
When viewed as functions of $c$ it is easy to see from
\eqref{quasi2} that both sides of \eqref{V1211}
satisfy $f(c)=f(cp^2)$. Consequently it is enough to give a
proof for $c=q^{n-m+1}$ with $m$ an integer
such that $m\geq 2n+1$. But this is nothing but
\eqref{new2} with $n\to m$ and $a\to bq^{-n-1}$.
\end{proof}

\subsection*{Acknowledgements}
I thank Vyacheslav Spiridonov for prompting me to
look for a proof of \eqref{biba} beyond induction,
and Michael Schlosser for helpful discussions.
I thank Mizan Rahman for pointing out \eqref{NR414}
and \eqref{RV78} in the literature.

\bibliographystyle{amsplain}

\begin{thebibliography}{99}

\bibitem{Andrews76}
G. E. Andrews,
\textit{On $q$-analogues of the Watson and Whipple summations},
SIAM J. Math. Anal. \textbf{7} (1976), 332--336.

\bibitem{AB02}
G. E. Andrews and A. Berkovich,
\textit{The WP-Bailey tree and its implications},
J. London Math. Soc. (2) \textbf{66} (2002), 529--549.

\bibitem{BW03}
A. Berkovich and S. O. Warnaar,
\textit{Positivity preserving transformations for $q$-binomial
coefficients},
arXiv:math.CO/0302320.

\bibitem{vDS00}
J. F. van Diejen and V. P. Spiridonov,
\textit{An elliptic Macdonald-Morris conjecture and multiple modular
hypergeometric sums},
Math. Res. Lett. \textbf{7} (2000), 729--746.

\bibitem{vDS01a}
J. F. van Diejen and V. P. Spiridonov,
\textit{Elliptic Selberg integrals},
Internat. Math. Res. Notices (2001), 1083--1110.

\bibitem{vDS01b}
J. F. van Diejen and V. P. Spiridonov,
\textit{Modular hypergeometric residue sums of elliptic Selberg integrals},
Lett. Math. Phys. \textbf{58} (2001), 223--238.

\bibitem{vDS02}
J. F. van Diejen and V. P. Spiridonov,
\textit{Elliptic beta integrals and modular hypergeometric sums: an overview},
Rocky Mountain J. Math. \textbf{32} (2002), 639--656.

\bibitem{FT97}
I. B. Frenkel and V. G. Turaev,
\textit{Elliptic solutions of the Yang-Baxter equation and modular
hypergeometric functions},
The Arnold-Gelfand mathematical seminars, 171--204,
(Birkh\"auser Boston, Boston, MA, 1997).

\bibitem{GR90}
G. Gasper and M. Rahman,
\textit{Basic Hypergeometric Series},
Encyclopedia of Mathematics and its Applications, Vol.~35,
(Cambridge University Press, Cambridge, 1990).

\bibitem{KN03}
Y. Kajihara and M. Noumi,
\textit{Multiple elliptic hypergeometric series --An approach from
the Cauchy determinant--},
arXiv:math.CA/0306219.

\bibitem{KNR03}
E. Koelink, Y. van Norden and H. Rosengren,
\textit{Elliptic U(2) quantum group and elliptic hypergeometric series},
arXiv:math.QA/0304189.

\bibitem{NR81}
B. Nassrallah and M. Rahman,
\textit{On the $q$-analogues of some transformations of nearly-poised
hypergeometric series},
Trans. Amer. Math. Soc \textbf{268} (1981), 211--229.

\bibitem{RV93}
M. Rahman and A. Verma,
\textit{Quadratic transformation formulas for basic hypergeometric series},
Trans. Amer. Math. Soc \textbf{335} (1993), 277--302.

\bibitem{Rosengren01}
H. Rosengren,
\textit{A proof of a multivariable elliptic summation formula
conjectured by Warnaar}, in \textit{$q$-Series with Applications
to Combinatorics, Number Theory, and Physics}, pp. 193--202,
B.~C.~Berndt and K.~Ono eds.,
Contemp. Math. Vol.~291 (AMS, Providence, 2001).

\bibitem{Rosengren02}
H. Rosengren,
\textit{Elliptic hypergeometric series on root systems},
arXiv:math.CA/0207046.

\bibitem{RS03}
H. Rosengren and M. Schlosser,
\textit{Summations and transformations for multiple basic and
elliptic hypergeometric series by determinant evaluations},
arXiv:math.CA/0304249.

\bibitem{Slater66}
L. J. Slater,
\textit{Generalized hypergeometric functions},
(Cambridge University Press, Cambridge, 1966).

\bibitem{Spiridonov01}
V. P. Spiridonov,
\textit{Elliptic beta integrals and special functions of hypergeometric type},
in \textit{Integrable structures of exactly solvable two-dimensional
models of quantum field theory}, pp. 305--313,
S.~Pakuliak and G.~von Gehlen eds.,
NATO Sci. Ser. II Math. Phys. Chem. Vol.~35,
(Kluwer Academic Publishers, Dordrecht, 2001).

\bibitem{Spiridonov02a}
V. P. Spiridonov,
\textit{Theta hypergeometric series},
in \textit{Combinatorics with Applications to Mathematical Physics},
pp. 307--327, V.~A.~Malyshev and A.~M.~Vershik, eds.,
(Kluwer Academic Publishers, Dordrecht, 2002).

\bibitem{Spiridonov02}
V. P. Spiridonov,
\textit{An elliptic incarnation of the Bailey chain},
Int. Math. Res. Notices. \textbf{37} (2002), 1945--1977.

\bibitem{Spiridonov03}
V. P. Spiridonov,
\textit{Theta hypergeometric integrals},
arXiv:math.CA/0303205.

\bibitem{SZ99}
V. Spiridonov and A. Zhedanov,
\textit{Classical biorthogonal rational functions on elliptic grids},
C. R. Math. Acad. Sci. Soc. R. Can. \textbf{22} (2000), 70--76.

\bibitem{W02}
S. O. Warnaar,
\textit{Summation and transformation formulas
for elliptic hypergeometric series},
Constr. Approx. \textbf{18} (2002), 479--502.

\bibitem{W03}
S. O. Warnaar,
\textit{Extensions of the well-poised and elliptic well-poised
Bailey lemma},
Indag. Math. (N.S.), to appear.

\end{thebibliography}

\end{document}